# Eighty Years of the Finite Element Method: Birth, Evolution, and Future

**Wing Kam Liu, Northwestern University**
**Shaofan Li, UC Berkeley**
**Harold S. Park, Boston University**

This year marks the eightieth anniversary of the invention of the finite element method (FEM). FEM has become the computational workhorse for engineering design analysis and scientific modeling of a wide range of physical processes, including material and structural mechanics, fluid flow and heat conduction, various biological processes for medical diagnosis and surgery planning, electromagnetics and semi-conductor circuit and chip design and analysis, additive manufacturing, i.e. virtually every conceivable problem that can be described by partial differential equations (PDEs). FEM has fundamentally revolutionized the way we do scientific modeling and engineering design, ranging from automobiles, aircraft, marine structures, bridges, highways, and high-rise buildings. Associated with the development of finite element methods has been the concurrent development of an engineering science discipline called computational mechanics, or computational science and engineering.

In this paper, we present a historical perspective on the developments of finite element methods mainly focusing on its applications and related developments in solid and structural mechanics, with limited discussions to other fields in which it has made significant impact, such as fluid mechanics, heat transfer, and fluid-structure interaction. To have a complete storyline, we divide the development of the finite element method into four time periods: I. (1941-1965) Early years of FEM; II. (1966-1991) Golden age of FEM; III. (1992-2017) Large scale, industrial

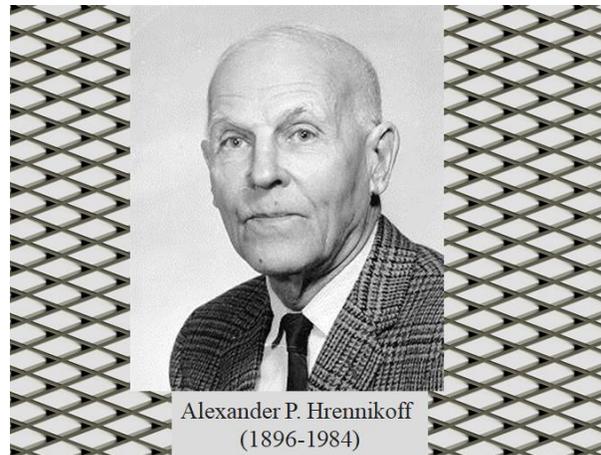

Alexander P. Hrennikoff
(1896-1984)

applications of FEM and development of material modeling, and IV (2018-) the state-of-the-art FEM technology for the current and future eras of FEM research. Note that this paper may not strictly follow the chronological order of FEM developments, because often time these developments were interwoven across different time periods.

## I.    (1941-1965) The birth of the finite element method

The origin of the finite element method as a numerical modeling paradigm may be traced back to in early 1940's. In 1941, A. Hrennikoff, a Russian-Canadian structural engineer at the University of British Columbia, published a paper in ASME *Journal of Applied Mechanics* on his membrane and plate model as a lattice framework. This paper is now generally regarded as a turning point that led to the birth of FEM. In Hrennikoff's 1941 paper, he discretized the solution domain into a mesh of lattice structure, which was the earliest form of a mesh discretization.

On May 3rd, 1941, the same year that Hrennikoff published his paper, R. Courant of New York University delivered an invited lecture at a meeting of the American Mathematical Society held in

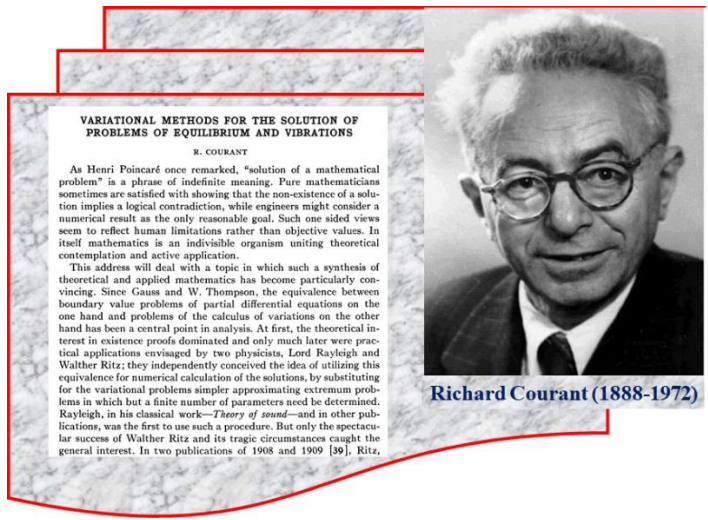

**Richard Courant (1888-1972)**

Washington D.C. on his numerical treatment using a variational method to solve a second order PDE, which arises from Saint-Venant's torsion problem of a cylinder. In this work, Courant systematically used the Rayleigh Ritz method with a trial function defined on finite triangle subdomains, which is a primitive form of the finite element method. Courant's presentation was later published as a paper in 1943. Similar works of discretization and variational formulations were also reported in the literature, including McHenry [1943], Prager and Synge [1947], and Synge [1957]. As Ray Clough commented in his 1980 paper, *One aspect of the FEM, mathematical modeling of continua by discrete elements, can be related to work done independently in the 1940s by McHenry and Hrennikoff-formulating bar element assemblages to simulate plane stress systems. Indeed, I spent the summer of 1952 at the Boeing Airplane Company trying to adapt this procedure to the analysis of a delta airplane wing, the problem which eventually led to the FEM.*

By the early 1950s, several engineers and academics had further developed and extended these early approaches to solve real engineering problems in aeronautical and civil engineering. In parallel but with different emphases, J. H. Argyris (1954) at the Imperial College London, and M. J. Turner (1950-1956) at Boeing Company, who was later joined by R. Clough of UC Berkeley and H.C. Martin of Washington University, developed what we know today as the earliest form of the finite element method (1954), which was called the Matrix Stiffness Method at the time. However, in 1960 R. Clough coined the phrase *Finite Element Method*, and this unassuming and right-to-the-point phrase was an instant hit, bringing out the essence of the method.

Using the language of engineers, Argyris "translated" Courant's variational approach into the energy method of engineering structures, while importantly, Turner, Clough, Martin, and Topp developed FEM interpolants for triangular elements, which is suitable for structural parts with arbitrary shape. In some sense, the invention of the triangle element was a "quantum leap", and hence for a large spectrum of the

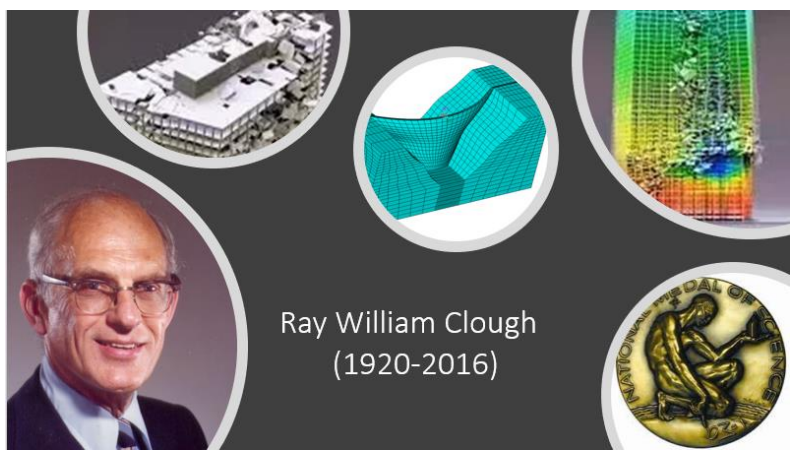

Ray William Clough (1920-2016)

engineering community, the inception of the FEM is the publication of the landmark paper by Turner, Clough, Martin, and Topp in 1956. The following is an excerpt from a 2014 document that celebrates the 50[th] anniversary of the formation of National Academy of Engineering, which is an official account of that part of FEM history: *To ensure safety and avoid costly modifications after planes entered production, engineers needed a reliable method for determine in advance whether their designs could withstand the stresses of fight. M. Jon Turner, head of Boeing's Structural Dynamics Unit, addressed that problem in the early 1950s by bringing civil engineering professor Ray Clough of the University of California, Berkeley, and Harold Martin of the University of Washington to Boeing for summer ``faculty internships,''*

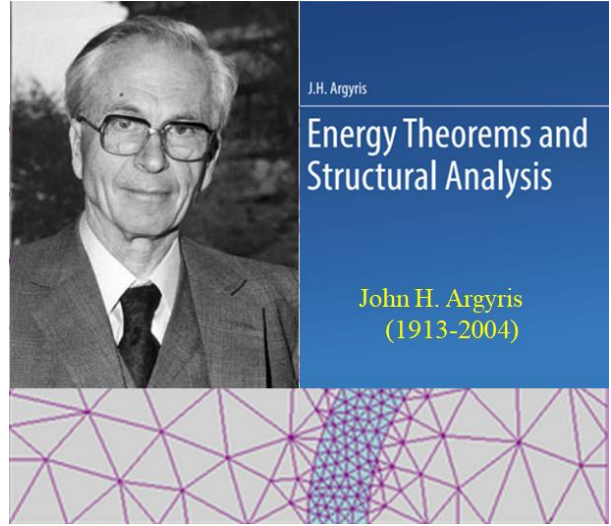

J.H. Argyris
**Energy Theorems and Structural Analysis**

John H. Argyris
(1913-2004)

*Collectively, they created a method of structural analysis that Turner applied at Boeing using computers to perform the myriad calculations need to predict real-world performance.* That fruitful collaboration led to Clough's development a few years later of what he named the finite element method (FEM). Clough formed a research group at UC Berkeley that used FEM in a host of analytical and experimental activities, from designing buildings and structures to withstand nuclear blasts or earthquakes to analyzing structural requirements for spacecraft and deep-water offshore drilling, led to Clough's development a few years later of what he named the finite element method (FEM). Clough formed a research group at UC Berkeley that used FEM in a host of analytical and experimental activities, from designing buildings and structures to withstand nuclear blasts or earthquakes to analyzing structural requirements for spacecraft and deep-water offshore drilling. By revolutionizing the application of computer technologies in engineering, FEM continues to help engineers design to this day all sorts of durable, cost-effective structures. Meanwhile, Turner's efforts at Boeing contributed to the success of its renowned line of commercial jets, beginning in 1958 with the 707 and continuing in 1964 with the 727, which could land on shorter runways and serve more airports. Equipped with three fuel-efficient turbofan engines, the 727 became the workhorse of commercial aviation and helped achieve a threefold increase in U.S. passenger air traffic in the '60s.

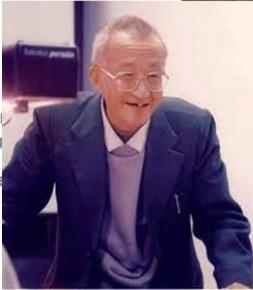

基于变分原理的差分格式

冯 康

Kang Feng
(1920-1993)

Independently and separately, in the early 1960s, Kang Feng of the *Chinese Academy of Science* also proposed a discretization numerical method based on variational principles for solving elliptic partial differential equations. As Peter Lax (1993) commented, "*Independently of parallel developments in the West. he (Feng) created a theory of the finite element method. He was instrumental in both the implementation of the method and the creation of its* *theoretical foundation using estimates in Sobolev spaces….*", which was one of the first convergence studies of finite element methods.

During this period, several great engineering minds were focusing on developing finite element methods. In particular, J. Argyris with his co-workers at the University of Stuttgart; R. Clough and colleagues such as E. L. Wilson and R.L. Taylor at the University of California, Berkeley; O. Zienkiewicz with his colleagues such as E. Hinton and B. Irons at Swansea University; P. Ciarlet at the University of Paris XI; R. Gallager and his group at Cornell University, R. Melosh at Philco Corporation, and B. Fraeijs de Veubeke at the University de Liege had made some important and significant contributions to developments of finite element methods.

To understand what happened sixty years ago, we quote an excerpt from a FEM history paper by Clough and Wilson in 1999, in which they recalled: *When Clough presented the first paper using the finite element terminology in 1960 it attracted the attention of his friend, Professor O. C. Zienkiewicz, who was then on the faculty at Northwestern University. A few weeks after the presentation of*

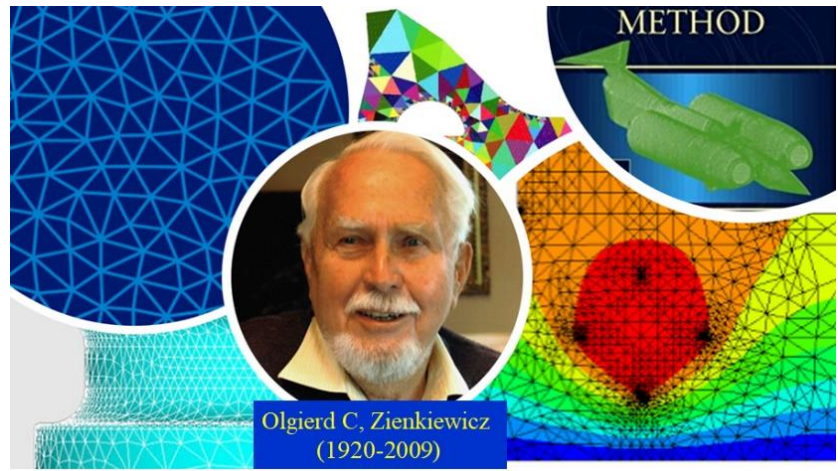
METHOD

Olgierd C. Zienkiewicz
(1920-2009)

*the paper Zienkiewicz invited Clough to present a seminar on the finite element method to his students. Zienkiewicz was considered one of the world's experts on the application of the finite difference method to the solution of continuum mechanics problems in Civil Engineering; therefore, Clough was prepared to debate the relative merits of the two methods. However, after a few penetrating questions about the finite element method, Zienkiewicz was almost an instant convert to the method.*

Zdenek Bazant, then a visiting Associate Research Engineer at UC Berkeley, recalled: …… *The founding of FEM was Clough's 1960 paper in the ASCE Conf. on Electronic Computation (Clough (1960), which was the <u>first</u> to derive, by virtual work, the finite element stiffness matrix of an element of a continuum (a triangular constant strain element). ……. I recall Ray Clough showing to me his 1962 report to the U.S. Engineer District, Little Rock, Corps of Engineers on his analysis of a crack observed in Norfolk Dam (Clough (1962)), during my stay in Berkeley in 1969. I was mesmerized by seeing that 1962 report. It presented a 2D stress analysis of large crack observed in Norfork dam. The dam was subdivided into about 200 triangular elements and provided stress contours for a number of loading cases. …… Clough was at that time way ahead of anybody else.*

To collaborate Z. Bazant's recollection, we cite Ed Wilson's recount of that part of the history:

*In 1956, Ray, Shirley, and three small children spent a year in Norway at the Ship Research Institute in Trondheim. The engineers at the institute were calculating stresses due to ship vibrations in order to predict fatigue failures at the stress concentrations. This is when Ray realized his element research should be called the Finite Element Method which could solve many*

*different types of problems in continuum mechanics. Ray realized the FEM was a direct competitor to the Finite Difference Method. At that time FDM was being used to solve many problems in continuum mechanics. His previous work at Boeing, the Direct Stiffness Method, was only used to calculate displacements not stresses.*

*In the fall semester of 1957, Ray returned from his sabbatical leave in Norway and immediately posted a page on the student bulletin board asking students to contact him if they were interested in conducting finite element research for the analysis of membrane, plate, shell, and solid structures. Although Ray did not have funding for finite element research, a few graduate students who had other sources of funds responded. At that time, the only digital computer in the College of Engineering was an IBM 701 that was produced in 1951 and was based on vacuum tube technology. The maximum number of linear equations that it could solve was 40. Consequently, when Ray presented his first FEM paper in September 1960, "The Finite Element Method in Plane Stress Analysis," at the ASCE 2nd Conference on Electronic Computation in Pittsburgh, Pennsylvania, the coarse-mesh stress-distribution obtained was not very accurate. Therefore, most of the attendees at the conference were not impressed. After the improvement of the speed and capacity of the computers on the Berkeley campus, Professor Clough's paper was a very fine mesh analysis of an existing concrete dam. The paper was first presented in September 1962 at a NATO conference in Lisbon, Portugal. Within a few months, the paper was republished in an international Bulletin, which had a very large circulation, as "Stress Analysis of a Gravity Dam by the Finite Element Method", (with E. Wilson), International Bulletin RILEM, No. 10, June 1963.*

*The Lisbon paper reported on the finite element analysis of the 250-foot-high Norfork Dam in Arkansas, which had developed a vertical crack during construction in 1942. The FEM analysis correctly predicted the location and size of the crack due to the temperature changes and produced realistic displacements and stresses within the dam and foundation for both gravity and several hydrostatic load conditions. Because of this publication, many international students and visiting scholars came to Berkeley to work with Professor Clough. Also, he freely gave the FORTRAN listing of their finite element analysis computer program to be used to evaluated displacement and stresses in other two-dimensional plane structures with different geometry, materials and loading. Therefore, professional engineers could easily use the powerful new FEM to solve for the stress distributions in their structural engineering problems in continuum mechanics. However, he did not capitalize on his success in the development of the FEM. He returned to the task of building the earthquake engineering program at Berkeley – the task he given when he was first hired in 1949.*

For his decisive contribution to the developments of FEM, Ray William Clough was awarded the National Medal of Science in 1994 by the then vice-president of the United States Al Gore. Today, the general consensus is that Ray Clough along with J. Argyris and O. C. Zienkiewicz made the most significant contributions to the developments of finite element method after an early mathematical pre-working of Richard Courant.

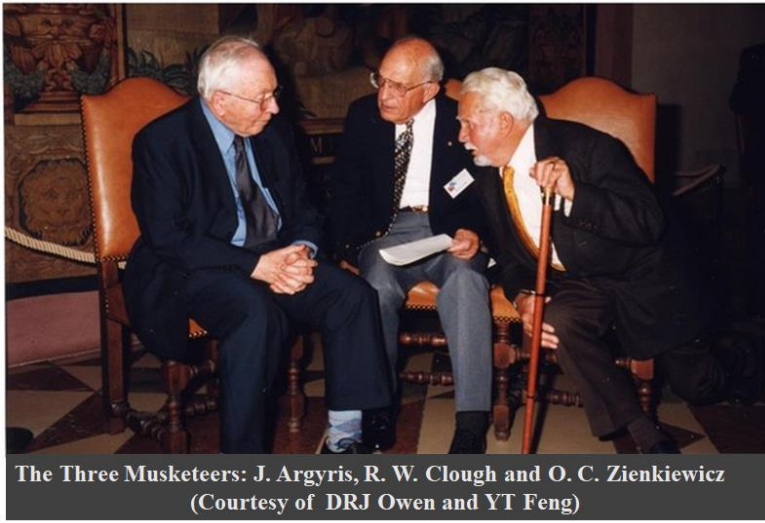

**The Three Musketeers: J. Argyris, R. W. Clough and O. C. Zienkiewicz**
**(Courtesy of DRJ Owen and YT Feng)**

It is worth noting that E.L. Wilson of UC Berkeley was the first person to develop finite element open-source software. An excerpt from the Clough and Wilson's paper in 1999 stated: *In 1958 Wilson, under the direction of Clough, initiated the development of an automated finite element program based on the rectangular plane stress finite element developed at Boeing. After several months of learning to program the IBM 701, Wilson produced a limited capacity, semiautomated program which was based on the force method. An MS research report was produced, which has long since been misplaced, with the approximate title of Computer Analysis of Plane Stress Structures. …… In 1959 the IBM 704 computer was installed on the Berkeley Campus. It had 32K of 32 bit memory and a floating point arithmetic unit which was approximately 100 times faster than the IBM 701. This made it possible to solve practical structures using fine meshes.*

It is also worth mentioning that in 1963, I.T. Oden, a senior structural engineering in the research and development division of General Dynamics Corporation at Fort Worth at the time, wrote a 163-page comprehensive technical report with G.C. Best, in which they developed a long list of solid and structure finite elements. including tetrahedral element, hexahedral element, thin plate element, thick plate element, plate element with stringers or stiffeners, composite sandwich plate elements, and shallow shell elements (Best and Oden (1963)). In fact, Oden and Best wrote one of the first general purpose finite element computer codes at the time. The Fortran FEM code developed by Oden and Best had an element library that includes elements for 3D elasticity, 2D plane elasticity, 3D beam and rod elements, composite layered plate and shell elements, and elements for general composite materials. Their work also included hybrid methods and stress based FEMs, which may be even earlier than those of Pian (Pian (1964)). Moreover, Oden and Best's FEM code is also able to handle FEM eigenvalue modal analysis, and numerical integrations over triangle and tetrahedra elements, and it had linear system solvers for general FEM static analysis that were among the most effective at that time. This FEM computer code was used for many years in aircraft analysis and design in the aerospace and defense industry.

In terms of worldwide research interest, by 1965, finite element research had become a highly active field, with the total number of papers published in the literature exceeding 1000. During this period, there were two seemingly unrelated events for FEM development, but significantly affected future FEM developments. These events were the discovery of mixed variational principles in elasticity. In 1950, E. Reissner rediscovered E. Hellinger's mixed variational principle from 1914, in which both the displacement field and the stress field are the primary unknowns. This variational principle is called the Hellinger-Reissner variational principle. Shortly after, H. Hu in 1954 and K. Washizu in 1955 proposed a three-field mixed variational principle in elasticity, which was called the Hu-Washizu variational principle. As early as 1964, Theodore H.H. Pian recognized the potential of using these variational principles to formulate Galerkin weak form-based finite element formulations and proposed the assumed stress finite element method. This began the use of mixed variational principles to formulate Galerkin finite element methods, which was followed by the assumed strain finite element method developed by Juan Simo and his co-workers in the later period of finite element developments.

## II.      (1966-1991) The Golden Age of the Finite Element Method

The mid 1960s saw rapid developments in finite element method research and applications. As TJR Hughes recalled, "*I first heard the words 'the Finite Element Method' in 1967 – which changed my life. I started to read everything that was available and convinced my boss to start the Finite Element Method Development Group, which he did. Dr. Henno Allik was Group Leader, and I was in the Group, then we added programmers. In one year, we had a 57,000-line code, GENSAM (Gen­eral Structural Analysis and Matrix Program, or something like that). That was 1969, and the code was continually developed thereafter and may still be in development and use at GD/Electric Boat and General Atomics, originally a division of General Dynamics*."

Starting from the end of 1960s, the rigorous approximation theory that underpins the finite element method started to be developed. This movement was first highlighted by the proof of optimal and superconvergence of finite element methods. This attracted the interests of some distinguished mathematicians all over the world, including G. Birkhoff, M.H. Schultz, R.S. Varaga, J. Bramble, M. Zlamal, J. Cea, JP. Aubin, J. Douglas, T. Dupont, LC. Goldstein, LR. Scott, J. Nitsche, AH. Schatz, PG. Cialet, G. Strang, G. Fix, JL. Lions, M. Crouzeix, PA. Raviart, and I. Babuska, AK. Aziz, and J.T. Oden. Some useful results developed for the proof of FEM convergence are the Cea lemma and the Bramble-Hubert lemma. The fundamental work of Nitsche [1970] on $L^\infty$ estimates for general classes of linear elliptic problems must stand out as one of the most important contributions for mathematical foundation of FEM in 1970s. It may be noted that unlike other mathematics movements, the convergence study of finite element methods was an engineering-oriented movement. The mathematicians soon found that, in practice, engineers were using either non-conforming FEM interpolants or numerical quadrature that violates variational principles or the standard bilinear form in Hilbert space. Gilbert Strang referred to these numerical techniques as "variational crimes". To circumvent complicated convergence proofs, the early patch tests were invented by Bruce Iron and Robert Melosh, which were proven to be instrumental for ensuring convergence to the correct solution.

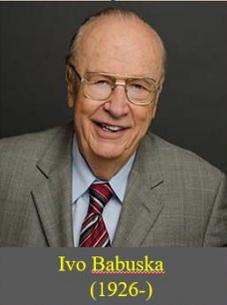

Following Theodore Pian's invention of the assumed stress element, attention shifted to the mixed variational principle based finite element method. In 1966, L. Hermann proposed a mixed variational principle for incompressible solids. However, most mixed variational principles are not extreme variational principles, and thus suffer from numerical instability. In early 1970s, I. Babuska and F. Brezzi discovered their groundbreaking results, known today as the Babuska-Brezzi condition, or the LBB condition, giving tribute to O. Ladyzhenskya --- a Russian female mathematician, who provided the early insight of this problem. The so-called LBB condition, or the inf-sup condition, provides a sufficient condition for a saddle point problem to have a unique solution that depends continuously on the input data; thus, it provides a guideline to construct shape functions for the mixed variational principle-based finite element method.

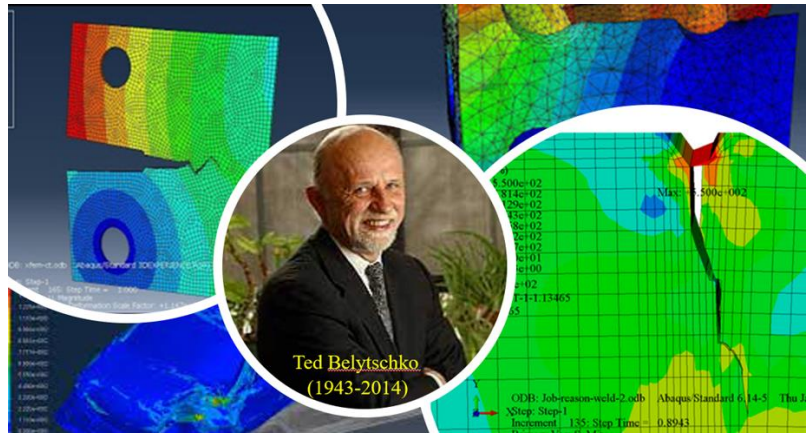

Entering the 1970s, FEM development began to focus on using FEM to simulate the dynamic behavior of structures, including crashworthiness in the automotive industry. Various time integration methods had been developed, including the Newmark-beta method, the Wilson-theta method, the Hilbert-Hughes-Taylor algorithm, the Houbolt integration algorithm, and the explicit time integration algorithms.

In the late 1970s, T. Belytschko, K. C. Park, and later TJR Hughes proposed using explicit or implicit-explicit, explicit-multiple explicit time integration, and implicit time integration with damping control to solve nonlinear structural deformation and structural dynamics problems. It turned out that the explicit time integration was a game changer for the automotive industry, establishing FEM technology as the main tool of passenger vehicle design and crashworthiness analysis. By the end of 1980s, there were thousands of workstations running explicit time integration-based FEM codes in the three major automakers in the United States.

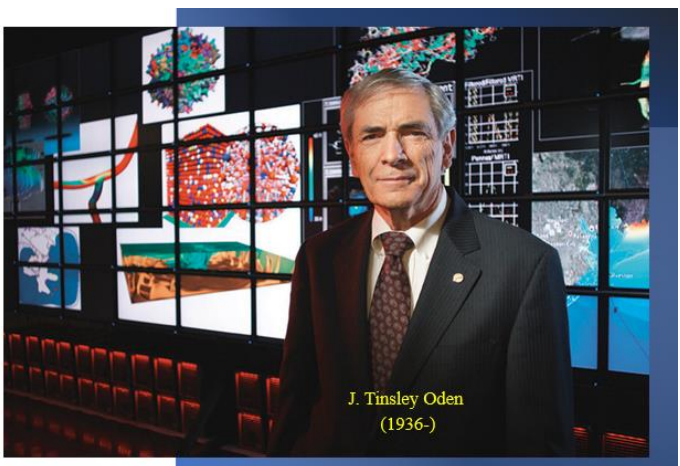

One of the main FEM research topics in the 1980s was using finite element techniques to solve the Navier-Stokes equation as an alternative to finite difference and finite volume methods. Starting the early 1970s, JT Oden and his co-workers had begun working on finite element solutions for fluid dynamics (Oden and Wellford (1972)). TJR Hughe and his co-workers Brooks and Tezduyar developed the streamline upwind/Petrov Galerkin method and later Stabilized Galerkin finite element method to solve Navier-Stokes equations under various initial and boundary conditions. Furthermore, Hughes and co-workers later developed space-time finite element methods and variational multiscale finite element methods (see Hughes and Tezdyuar (1981) Brooks and Hughes (1982), Mizukami and Hughes (1985), and Hughes (1995)). For this work, from 1986 to 1991, Hughes and his co-workers such as LP Franca and others wrote a ten-part series on finite element formulation for computational fluid dynamics (Hughes et al. (1986) and Shakib et al. (1991)).

Among the many advances in finite element technologies in 1980s, the most notable may belong to J. Simo at the University of California, Berkeley and later at Stanford University. Simo and Taylor developed the consistent tangent operator for computational plasticity, which was a milestone after the original concept of consistent linearization proposed by TJR Hughes and K. Pister. Moreover, after the Hughes-Liu 3D degenerated continuum shell and beam elements and the Belytchko-Tsay single-point element, Simo and his co-workers, such as L. Vu-Quoc and D.D. Fox, developed geometrically exact beam and shell theories and its finite element formulation. Simo and his co-worker such as MS Rifai and F. Armero also developed various assumed strain or enhanced strain methods for mixed variational formulations. It should be noted that E. Ramm and his colleagues at the University of Stuttgart have also made significant contributions on geometrically nonlinear shell element formulations over a span of more than thirty years e.g., Bischoff and Ramm (1998).

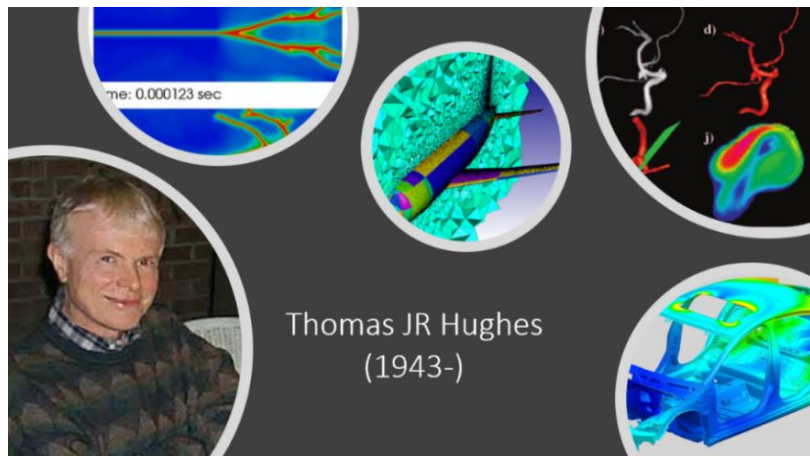

Another highlight of FEM technology is the development of FEM solvers for fluid-structure interaction. From the mid-1970s to early 1990s, there was an urgent need to develop techniques to solve large-scale fluid-structure interaction problems in the aerospace and civil engineering industries. A class of finite element fluid-structure interaction solvers were developed, and some early contributors include J. Donea, A. Huerta (see Donea and Huerta (2003), Donea et al. (2017)), and the Hughes-Liu-Zimmermann Arbitrary Lagrangian-Eulerian (ALE) fluid-structure finite element formulation, which describes the moving boundary problem. ALE-based finite element simulations were used due to their ability to alleviate many of the drawbacks of traditional Lagrangian-based and Eulerian-based finite element formulations.

When using the ALE technique in engineering modeling and simulations, the computational mesh inside the domains can move arbitrarily to optimize the shapes of elements, while the mesh on the boundaries and interfaces of the domains can move along with materials to precisely track the boundaries and interfaces of a multi-material system. The invention of ALE finite element method may be credited to Hirt, Amsden, and Cook (1974). C. Farhat was the first person to use large scale parallel ALE-FEM solver to compute fluid-structure interaction problems. He and his group systematically applied FEM-based computational fluid dynamics (CFD) solvers for aircraft structure design and analysis. They developed the finite element tearing and interconnecting (FETI) method for the scalable solution of large-scale systems of equations on massively parallel processors.

The finite element fluid-structure interaction research had a major impact on many practical applications, such as providing the foundation for the patient specific modeling of vascular disease and the FEM-based predictive medicine later developed by CA Taylor and TJR Hughes and their co-workers in the mid-1990s. Holzapfel, Eberlein, Wriggers, and Weizsäcker developed large strain finite element formulation for soft biological membranes (Holzapfel et al. (1996)).

Another major milestone in the development of finite element methods was the invention and the development of nonlinear probabilistic or random field finite element method, which was first developed by W. K. Liu and T. Belytschko in the late 1980s. By considering uncertainty in loading conditions, material behavior, geometric configuration, and support or boundary conditions, the probabilistic finite element method provided a stochastic approach in computational mechanics to account for all these uncertain aspects, which could then be applied in structure reliability analysis. The random field finite element research has become crucial in civil and aerospace engineering and the field of uncertainty quantification.

In the early 1980's, M.E. Botkin at General Motors research Lab, and N. Kikuchi and his group at the University of Michigan developed structural shape optimization finite element method for the automotive industry. Other contributors include M. H. Imam from Uman Al-Qura University. This preceded the seminal 1988 paper of Bendsoe and Kikuchi, who developed a homogenization approach to finding the optimal shape of a structure under prescribed loading. Later developments in topology optimization were driven by O. Sigmund and his students, such as Gengdong Cheng in China.

One of the driving forces in the finite element method development during in the early decades was the safety analysis of big dams, at first, then of concrete nuclear reactor vessels for gas-cooled reactors, e.g., J. Rashid (1967), nuclear containments, and hypothetical nuclear reactor accidents (Marchertas et al. (1978), Bazant et al. (1980)) and of tunnels and of foundations for reinforced concrete structures. To simulate concrete failure, the vertical stress drop in finite elements and progressive softening technique were introduced already in 1968. However, the spurious mesh sensitivity and the necessity of localization limiter was generally overlooked until demonstrated mathematically in Bazant's (1976), based on the fact that strain softening states of small enough test specimens in stiff enough testing frames are stable. Numerically this was demonstrated by crack band finite element calculations in Bazant and Cedolin (1978), where it was also shown that spurious mesh sensitivity of sudden stress drop can be avoided by adjusting the material strength so as to ensure correct energy release rate. Hillerborg et al. (1976) avoided mesh sensitivity by using an interelement cohesive softening called the fictitious crack model. Despite the success of the early FE calculations, the concept of progressive strain-softening damage was not generally accepted by mechanicians until its validity and limitations were demonstrated by Bazant and Belytschko (1985)). They showed that the existence of elastic unloading stiffness (previously ignored) makes waves propagation in a strain-softening state possible. In 1989, Lubliner, Oliver, Oller, and Onate (Lublinear et al. (1989)) developed a plastic-damage theory-based finite element formulation to model concrete materials by introducing internal variables, which has the capacity of modeling concrete material degradation and cracking. Today, the multiscale based homogenization and damage analysis method is the state-of-the-art finite element modeling for concrete materials (Wriggers and Moftah (2006)).

In the middle of 1980s, the mesh sensitivity issue in calculating strain softening or strain localization problem in computational plasticity became a challenging topic. It was eventually accepted that the partial differential equations that are associated with the classical plasticity become ill-posed after the material passes the yield point and enters into the softening stage. This especially became a dire situation when civil engineers applied FEM to solve complex structural and geotechnical engineering problems, which involved complex plastic deformations of concrete, rock, soil, clay, and granular material in general. Bazant and others realized that this is because the classical continuum plasticity theory lacks an internal length scale. To remedy this problem, starting from the middle 1980s. Many efforts were devoted to establishing finite element formulations of nonlocal (Bazant, Belytschko, Chang (1984)), strain-gradient, strain-Laplacian media, micropolar or Cosserat continua, because they provided an internal length scale, allowing FEM simulations to capture, in a mesh-independent manner, the strain softening, strain localization or shear band formation. Pijaudier-Cabot and Bazant (1987) developed an effective

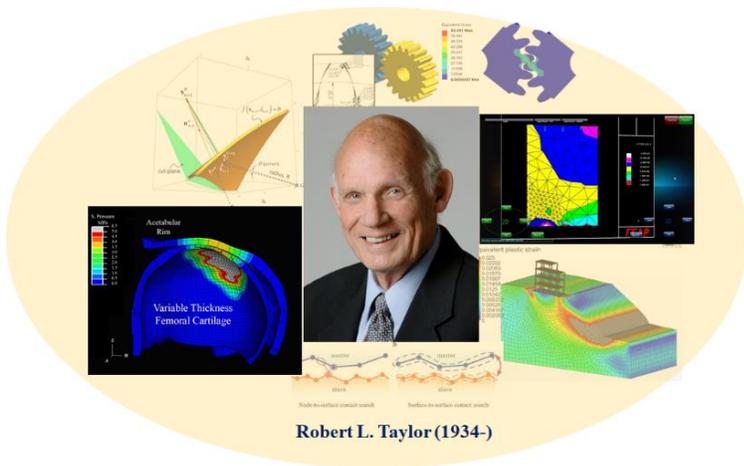

**Robert L. Taylor (1934-)**

nonlocal finite element method in which nonlocality is applied only to the damage strain.Other influential and representative works in this topic are from R. De Borst [1991,1992], Peerlings and de Borst (1999), K. Willam, and P. Steinmann, and their co-workers [1993,1998].

In 1976, TJR Hughes, RL. Taylor, JL. Sackman, A. Curnier, and W. Kanoknukulchai published a paper entitled "A finite element method for a class of contact-impact problems." (Hughes et al. (1976)). This is one of the earliest finite element analysis in computational contact mechanics. It is the very first work on FEM modeling of dynamic contact and impact problems, and it plays an important role in the simulation accuracy for engineering problems involving interaction between different continuum objects. Examples include sheet metal forming, target impact and penetration, and interaction between pavement and tires. Developing accurate finite element contact algorithms has been a focal point since the 1980s. Various FEM contact algorithms have been developed, and some main contributors are N. Kikuchi, T. Oden, J. Simo, P, Wriggers, R. Taylor, P. Papadopoulos, and TA. Laursen. FEM contact algorithm research remained an active research topic until late 1990s (see Kikuchi and Oden (1988), Simo, Wriggers, and Taylor (1985), Papadopoulos and Taylor(1992).

In this period, one interesting emerging area was the development of finite element exterior calculus by D. Arnold, R. Falk, and R. Winther [2006]. The FEM exterior calculus uses tools from differential geometry, algebraic topology, and homological algebra to develop FEM discretization that are compatible with the underlying geometric, topological, and algebraic structures of the problems that are under consideration.

## III.    (1992-2017) Broad Industrial Applications and Materials Modeling

The first major event in FEM development this period was the formulation of the Zienkiewicz-Zhu error estimator [1992], which was a major contribution to the mathematical approximation theory of finite element methods in the 1990s. The Zienkiewicz-Zhu posteriori error estimators provide the quality control of a finite element solution with an optimal use of computational resources by refining the mesh adaptively. The idea and spirit of Zienkiewicz-Zhu was further carried out by Ainsworth and Oden [1997][2000], and today using posteriori error estimation to improve the quality has been elevated to the height of Bayesian inference and Bayesian update. This research topic is now intimately related with what is now called validation and verification (V&V) procedure.

Since the late 1970s, Szabo [1978] and Babuska [1981,1982] started to develop *hp* version of FEM based on piecewise-polynomial approximations that employ elements of variable size *(h)* and polynomial degree *(p)*. They discovered that the finite element method converges *exponentially* when the mesh is refined using a suitable combination of h-refinements (dividing elements into smaller ones) and p-refinements (increasing their polynomial degree). The exponential convergence makes the method a very attractive choice compared to most other finite element methods, which only converge with an algebraic rate. This work continued until the late 1990s, spearheaded by M. Ainsworth, L. Demkowicz, JT. Oden, CAM. Duarte, OC. Zienkiewicz, and CE. Baumann (see Demkowicz et al. (1989) Oden et al. (1989)). J. Fish [1992] also proposed a s-version FEM by superposing additional mesh(es) of higher-order hierarchical elements on top of the original mesh of $C^0$ FEM discretization, so that it increases the resolution of the FEM solution.

To solve material and structural failure problems, research work in the 1990s focused on variational principle based discretized methods to solve fracture mechanics problems or strain localization problems. In 1994, Xu and Needleman developed a finite element cohesive zone model (CZM) that can simulate crack growth without remeshing, which was later further improved by M. Ortiz and his co-workers, who used CZM finite element method to solve fragmentation and material fatigue problems. It should be noted that long before the invention of cohesive zone model, S.T. Pietruszczak and Z. Mroz (1981) developed the first cohesive finite element for shear fracture in soil. Later, Bazant's group at Northwestern University developed various interface finite element methods, such as the microplane model, to study size effects of concrete material and other composite materials (Brocca and Bazant (2001); Caner and Bazant (2009)). Such models became a standard tool for simulating missile impact and explosions at, e.g., ERDA Vicksburg. A microplane FE model for fiber composites has been developed for Chrysler and Ford Co. to compare various designs of automobile crush-cans (Smilauer et al. (2011)). An anisotropic poromechanical microplane model has been formulated and used for FE analysis of hydraulic fracturing (Rahimi et al. (2019)).

To alleviate mesh bias issues in modeling material fracture and damage problems, T. Belytschko and WK. Liu developed meshfree particle methods, namely the element-free Galerkin (EFG) method and the reproducing kernel particle method (RKPM), which are based on the moving least square method and the wavelet multiresolution analysis, respectively. RKPM provides consistency and thus convergence enhancements as compared to the popular smoothed particle hydrodynamics (SPH) method. S. Li [1999] and others have successfully employed meshfree Galerkin methods to accurately simulate adiabatic shear band formation and propagation with minimum mesh adaptation. At the same time, CA. Duarte and T. Oden developed the so-called hp-Cloud method; E. Onate, S. Idelsohn. OC. Zienkiewicz, and RL. Taylor developed a finite point method [1996], and SN. Atluri and T. Zhu [1998] proposed a meshless local Petrov-Galerkin (MLPG) method, among many other meshfree methods. Fleming and Belytschko also showed that singularity functions could be included in the approximation functions to greatly improve simulations involving fracture mechanics (Fleming et al. (1997)).

One of the most challenging problems in the development of meshfree Galerkin method is how to integrate the weak form, because the meshfree interpolants are highly irregular and it is difficult to make them variational consistent. In 2001, J. S. Chen and his co-workers proposed a stabilized conforming nodal integration method for meshfree RKPM method, which is not only simple and stable, but also variational consistent with the Galerkin weak formulation (Chen et al. (2001)).

Shortly after the meshfree method developments, I. Babuska and his co-workers developed the partition of unity finite element method (PUFEM), which was later coined generalized finite element method (GFEM) (see Melenk and Babuska (1996) and Babuska and Melenk (1997)). PUFEM is a powerful method because it can be used to construct FE spaces of any given regularity, which is a generalization of the h, p, and hp versions of the finite element method, as well as providing the ability to embed an analytic solution into the FEM discretization instead of relying upon a generic polynomial interpolant.

A significant breakthrough in computational fracture mechanics and finite element refinement technology came in the late 1990s, when Ted Belytschko and his co-workers, including T. Black, N. Moes, and J. Dolbow, developed the eXtended finite element (X-FEM) (see Belytschko and Black (1999), and Moës et al. (1999)), which uses various enriched discontinuous shape functions to accurately capture the morphology of a cracked body without remeshing. Because the adaptive enrichment process is governed by the crack tip energy release rate, X-FEM provides an accurate solution for linear elastic fracture mechanics (LEFM). In inventing X-FEM, Ted Belytschko brilliantly utilized the concept PUFEM in solving fracture mechanics problems without remeshing. Entering the new millennium, A. Karma and his co-workers (Karma et al. [2000]) proposed the phase-field finite element method to solve crack growth and crack propagation problems, as the phase field method can accurately predict material damage for brittle fracture. Some of the leading contributors for this research topic are B. Bourdin, TJR. Hughes, C. Kuhn, R. Muller, C. Miehe, CM. Landis, among others.

As mentioned before, the main reason for the huge success of finite element methods is their broad applicability to engineering analysis and design across scientific disciplines. On the other hand, most mechanical engineering designs are performed by using various computer-aided design (CAD) tools, such as solid modeling. To directly blend the finite element method into CAD design

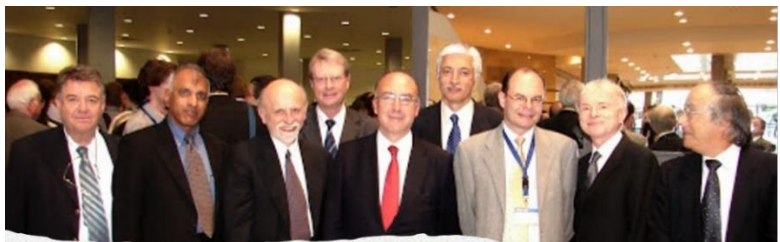

Gathering of International association for computational mechanics tools, TJR. Hughes and his co-workers such as JA Cottrell and Y. Bazilevs (Hughes et al. (2005), Cottrell et al. (2009)) developed the isogeometric analysis (IGA) finite element method, which established the Galerkin variational weak formulation in the control mesh and uses the Non-uniform rational basis spline (NURBS) functions as the FEM shape function to solve the problem at design stage.

Due to the emergence of nanotechnology, various multiscale methods have been developed to couple atomistic methods such as molecular dynamics and density functional theory (DFT) and other ab initio methods with continuum scale finite element methods. The most notable contributions in this area are hand-shake method (F.F. Abraham, J. Q. Broughton, N. Bernstein, and E. Kaxiras (1998)) quasi-continuum finite element method (E. B. Tadmor, M. Ortiz, R. Phillips (1996)), and the bridging scale method developed by G. Wagner and WK. Liu (2003), which is an ingenious use of the Mori-Zwanzig formalism in multiscale coupling technique. In 2007, V. Gavini, K. Bhattacharya, and M. Ortiz developed quasi-continuum orbital-free density-functional theory (DFT) finite element formulation for multi-million atom DFT calculations. (Gavini et al. (2007)).

Relatedly, there has been a rapid development in computational homogenization methods since the late 1990s, which is an alternate approach to obtaining continuum-scale properties based on smaller scale microstructures. The computational homogenization method or the finite element homogenization methods for composite materials may be divided into two main categories: (1) Computational asymptotic homogenization method or multiscale computational homogenization, which is aimed for modeling composite materials with periodic microstructure. The pioneer of multiscale homogenization finite element method may be credited to TY Hou and his co-workers (1997), and other earlier contributors are N. Kikuchi, S. Ghosh, J. Fish and M. S. Shephard; and (2) Computational micromechanics method, which is mainly aimed for composite materials with random microstructure, even though it may also be applied to materials with periodic microstructure, where the main contributors are: P. Suquet and his colleagues at the French National Centre for Scientific Research (CNRS), G. J. Dvorak (1994) and J. Fish at Rensselaer Polytechnic Institute and C. Miehe (1999) and his group at the University of Stuttgart. During the same period, many novel FEM fluid-structure interaction solvers have been developed, for instance, the immersed finite element method developed by L. Zhang and her co-workers (2004), which was motivated by the immersed boundary method pioneered by C. Peskin of Courant Institute of Mathematical Sciences, New York University.

An important finite element application area emerged with the developments of computational plasticity. The early finite element computational plasticity formulation is based on hypoelastic-plastic rate formulation. To satisfy the objectivity requirement, TJR Hughes and J. Winget first proposed the so-called incremental objectivity (1980), which was in the first time the esoteric continuum mechanics theory was applied to finite element formulations and computations, and it in turn promoted the development of nonlinear continuum mechanics in 1980s and 1990s. Soon afterwards. Simo and Hughes then extended Hughes-Winget incremental objectivity algorithm to finite deformation case in computational plasticity. The notion of consistency between the tangent stiffness matrix and the integration algorithm employed in the solution of the incremental problem was introduced by Nagtegaal (1982) and Simo and Taylor (1985). Consistent formulations have been subsequently developed for finite deformation plasticity by Simo and Oritz (1985,1988) within the framework of multiplicative decomposition of the deformation gradient and hyperelasticity. Also in 1980s, based on the Gurson's model, Tvergaard and Needleman (1984) developed the finite element based Gurson-Tvergaard-Needleman model, which is probably the most widely used finite element computational plasticity constitutive model used in material modeling, though recently attention has turned to developing machine-learning based or data-driven computational plasticity models e.g., E. Chinesta, et al (2017), and the unsupervised machine learning data-driven finite element methods, call Self-consistent Clustering Analysis (SCA), by Z. Liu et al. (2016, 2018).

An important advance of finite element method is the development of crystal plasticity finite element method (CPFEM), which was first introduced in a landmark paper by D. Piece, RJ Asaro, and A. Needleman in 1982. In the past almost four decades, there are a numbers of people who have significant contributions to the subject, for example, A. Arsenlis and DM. Parks from MIT and Lawrence Livermore National Laboratory (1999,2002,2004); PR. Dawson and his colleagues at Cornell University (1989,2011), and D. Raabe and his colleagues at Max-Planck-Institute fur Eisenforschung (2003,2004,2006), among others. Based on crystal slip, CPFEM can calculate dislocation, crystal orientation and other texture information to consider crystal anisotropy during computations, and it has been applied to simulate crystal plasticity deformation, surface roughness, fractures and so on. Recently, S. Li and his co-workers developed a finite element based multiscale dislocation pattern dynamics to model crystal plasticity in single crystal (2019). Yu et al., (2019) reformulated the self-consistent clustering analysis (SCA) for general elasto-viscoplastic materials under finite deformation. The accuracy and efficiency for predicting overall mechanical response of polycrystalline materials is demonstrated with a comparison to traditional full-field finite element methods.

In 2013, a group of Italian scientists and engineers led by L. Beirão da Veiga and F. Brezzi proposed a so-called virtual element method (VEM). The virtual element method is an extension of the conventional finite element method for arbitrary element geometries. It allows polytopal discretizations (polygons in 2-D or polyhedra in 3-D), which may be even highly irregular and non-convex element domains. The name *virtual* derives from the fact that knowledge of the local shape function basis is not required, and it is in fact never explicitly calculated. VEM possesses features that make it superior to the conventional FEM for some special problems such as the problems with complex geometries for which a good quality mesh is difficult to obtain, solutions that require very local refinements, and among others. In these special cases, VEM demonstrates robustness and accuracy in numerical calculations, when the mesh is distorted.

As early as 1957, R. Clough introduced the first graduate finite element course in UC-Berkeley, and since then finite element courses at both graduate and undergraduate levels have been added into engineering higher education curriculums in all the major engineering schools and Universities all over the world. As JT Oden recalled in his 1963 paper, ``*I went on to return to academia in 1964 and among my first chores was to develop a graduate course on finite element methods. At the same time, I taught mathematics and continuum mechanics, and it became clear to me that finite elements and digital computing offered hope of transforming nonlinear continuum mechanics from a qualitative and academic subject into something useful in modern scientific computing and engineering.*''
By the end of 2015, there have been more than several hundred of finite element methods monographs and textbooks published in dozens of languages worldwide. An exposition of finite element mathematical theory by Strang and Fix [1972] was among the earliest of those FEM books. JT Oden and his collaborators such as JN Reddy and GF Carey wrote a five-volume finite element monographs in late 1970s and early 1980s. Other FEM books such as those by Zienkiewicz and Cheung, or later Zienkiewicz and Taylor; Cook, Malkus, Plesha, and Witt; TJR. Hughes; KJ. Bathe, and the Nonlinear FEM monographs by Belytschko, Liu, and Moran (1999) and by De Borst,

Crisfield, Remmers, and Verhoosel (2012), among others, all have made major impacts on FEM educations and applications. Among all these FEM monographs and textbooks, the book by Zienkiewicz and Taylor or Zienkiewicz, Taylor, and Zhu probably has had most impacts on FEM technology popularization, which may be because R. L. Taylor wrote a FEM research computer program code named FEAP, which was placed in the appendix of that book, providing an immediate guidance on hands-on experience of implementing FEM for the readers.

The development of finite element software technology started in the 1960s. In 1963, E.L. Wilson and R. Clough developed a structural mechanics FEM code called Symbolic Matrix Interpretive System, SMIS, and then Wilson initiated and developed a general-purpose static and dynamic Structural Analysis Program, SAP. In late 1960s and the early 1970s, KJ. Bathe developed nonlinear FEM code ADINA based on SAP IV and NONSAP. Today, the brand name SAP2000 has become synonymous with the state-of-the-art finite element structural analysis and design methods since its introduction over 55 years ago. At the same period, NASA developed its own FEM code called NASTRAN. About the same time in late 1960s, J. Swanson asked his employer Westinghouse to develop FEM computer code, and his suggestion was rejected, and then he left the company and developed the initial ANSYS FEM code. Several years later, J.O. Hallquist at Lawrence Livermore National Laboratory also developed a 3D nonlinear FEM code called DYNA3D, which later evolved to LS-DYNA. By the end of 1990s and early 2000s, and FEM software industry has become a multi-billion-dollar business. There were several household finite element software company names such as ANSYS, ABAQUS, ADINA, LS-DYNA, NASTRAN, COMSOL Multiphysics, CSI, among others. Today, there are also a plethora of open-source FEM software available online, such as FreeFEM, OpenSees, Elmer, FEBio, FEniCS Project, DUNE, among some others.

## IV. (2018- present) Coming of A New Era

The modern form of the finite element method can routinely solve many industrial problems. They enable fundamental understanding and allow predictive analysis for product design. For new scientific discovery and engineering innovation, the development of new scientific principles often trails the pace of new inventions with the sheer volume of data that are generated across multiple spatial, temporal, and design parameters (spatial-temporal-parameter) spaces. For this reason, FE researchers are studying various forms of machine and deep learning methods, of which this class of methods cover the largest class of interpolations. According to the universal approximation theorem, a neural network (NN) can be designed and trained to approximate any given continuous function to desired accuracy [Selmic and Lewis (2002), Hashash et al. (2004)] which is believed to drive new discoveries and enable future computational discretization technologies. In this context, Mechanistic Data Science (MDS) finite element methods, which combine known scientific principles with newly collected data, will provide the critically needed research that can be a boon for new inventions.

Scientific and engineering problems typically fall under three categories: (1) problems with abundant data but undeveloped or unavailable scientific principles, (2) problems that have limited data and limited scientific knowledge, and (3) problems that have known scientific principles with uncertain parameters, with possible high computational load (Saha et al., 2021). In essence, mechanistic data science (MDS) mimics the way human civilization has discovered solutions to difficult and unsolvable problems from the beginning of time. Instead of heuristics, MDS uses machine learning methods like active deep learning and hierarchical neural network(s) to process input data, extract mechanistic features, reduce dimensions, learn hidden relationships through regression and classification, and provide a knowledge database. The resulting reduced order form can be utilized for design and optimization of new scientific and engineering systems (Liu et. al., 2021). Thus, the new focus of the finite element research has shifted towards the development of machine learning based finite element methods and reduced order models.

With the recent development of machine learning and deep learning methods, solving FEM by constructing a deep neural network has become a state-of-the-art technology. Earlier research focused on building up a shallow neural network following the FEM structure to solve boundary problems. J. Takeuchi and Y. Kosugi in 1993 proposed a neural network representation of the finite

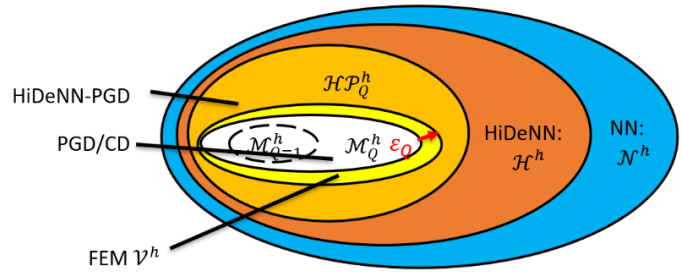

Approximation function spaces of FEM, PGD/CD, HiDeNN-PGD, HiDeNN, and NN, by Zhang et al. (2021)

element method to solve Poisson equation problems. G. Yagawa and O. Aoki (1995) replaced the FEM functional with the network energy of interconnected neural networks (NNs) to solve a heat conduction problem. Due to the limitation of computational power and slow convergency rate in shallow neural networks, earlier applications could only solve simple PDE problems. After the 2010s, neural networks for solving computational mechanics problems has become increasingly popular the rapid growth of deep learning techniques and the development of more sophisticated neural network structures, such as convolutional neural networks (CNN), Generative Adversarial Networks (GAN), residual neural networks (ResNet). For its high dimensional regression ability, some researchers, for example, F. Ghavamiana and A. Simone (2019), use deep neural network as a regression model to learn the material behavior or microstructure response. Other works focus on solving PDEs using deep learning neural networks. G. Karniadakis and his coworkers (2017, 2019, 2021) proposed a Physics-Informed Neural Networks (PINNs) to solve high dimensional PDEs in strong form that with constraints to accommodate both natural and essential boundary conditions. The idea of constructing deep neural network following the FEM structure is investigated again with advanced neural network methodology. Weinan E and B. Yu (2017) proposed a Deep Ritz Method for solving variational problems. J. Sirignano and K. Spiliopoulos (2018) proposed the so-called Deep Galerkin Method (DGM) to solve high-dimensional PDEs. Zabaras and his coworkers (2019) proposed a CNN-based physics-constrained deep learning framework for high-dimensional surrogate modeling and uncertainty quantification. T. Rabczuk and his coworkers (2020) systematically explore the potential to use NNs for computational mechanics by solving energetic format of the PDE. N. Trask and his coworkers (2021) proposed a partition of unity network for deep hp approximation of PDEs and extensively the training and initialization strategy to accelerate the convergence of the solution process. The constructing of element shape function by activation functions has been studied by J. Opschoor and his coworkers (2017) and J. He and his coworkers (2020).

Inspired by the universal approximation of deep neural networks (DNN), Zhang et al. (2020) published the first paper on the construction of the conventional finite element shape functions based in the hierarchical nature of the DNN, called **Hi**erarchical **De**ep-learning **N**eural **N**etworks (HiDeNN). Specifically, the authors demonstrated the construction of a few classes of deep learning interpolation functions, such as the reproducing kernel particle method (RKPM), non-uniform rational B-spline (NURBS), and isogeometric analysis (IGA), among other approximation techniques. Saha et al. (2021) generalized HiDeNN to a unified Artificial intelligence (AI)-framework, called HiDeNN-AI. HiDeNN-AI can assimilate many data-driven tools in an appropriate way, which provides a general approach to solve challenging science and engineering problems with little or no available physics as well as with extreme computational demand.

To reduce the FE computational cost, the so-called two-stage data-driven methods have been proposed of which during the offline stage, a data base generated by the finite element methods is first developed, and the final solutions are computed during the online stage. C. Farhat and his group at Stanford University have developed several dimensional reduction of nonlinear finite element dynamic models, including mesh sampling and weighting for the hyper reduction of nonlinear Petrov-Galerkin reduced-order models. To further reduce the finite element computational burden, Z. Liu et. al., (2016) applied the unsupervised machine learning techniques, such as the k-mean clustering method to group the material points during the offline stage and obtain the final solutions by solving the reduced-ordered Lippmann-Schwinger micromechanics equations. This class of data-driven approaches circumvents the computational burden of the well-established FE square method, the offline-online database approach to solve the concurrent FEM problems. They named the method the Finite Element-self-consistent clustering analysis (FE-SCA) of which the computational cost of the microscale analysis is reduced tremendously in multiple orders of magnitude speed up. J. Gao proposed an alternative (FE-SCAxSCA…xSCA) clustering analysis, of which the continuum finite element scale is concurrently solved with the (n-1) coupled-scale Lippmann-Schwinger micromechanics equations (Gao, 2020).

M. Ortiz and his group at Caltech developed data-driven finite elements for dynamics and noise data. J.-S. Chen and his co-workers have developed a physics-constrained data-driven RKPM method based on locally convex reconstruction for noisy databases (He and Chen 2020). S. Li and his group at UC-Berkeley utilized finite element solution generated data to develop a machine learning based inverse solution to predict pre-crash data of car collision. Bessa et al., (2017) proposed a data-driven framework to address the longstanding challenge of a two-scale analysis and design of materials under uncertainty applicable to problems that involve unacceptable computational expense when solved by standard FEM analysis of representative volume elements. The paper defined a framework that incorporates the SCA method to build large databases suitable for machine learning. The authors believe that this will open new avenues to finding innovative materials with new capabilities in an era of high-throughput computing ("big-data").

Reduced order modeling has been another active research field over the last decades. Early research works focused on the proper orthogonal decomposition (POD) method (also known as Karhunen-Loève transform, or principal component analysis) with the purpose of reducing the degrees of freedom of the discretized equations. The POD based model reduction has shown great success in computational fluid dynamics, see e.g., the works of P. Holmes and J. L. Lumley (1993). For further accelerating the simulations, K. Willcox and her coworkers (2008, 2002) proposed a missing point estimation method, which is known later as a hyper reduction method. Other notable works related to POD and hyper reduction methods are the Gauss–Newton with approximated tensors (GNAT) method, Grassmann manifold based reduced basis adaptation, thanks to C. Farhat and his coworkers (2008, 2013). For solid mechanics, D. Ryckelynck and his group (2005) proposed a hyper reduction method based on FEM for dealing with nonlinear problems. Another type of model reduction method, which is based on mathematics and has a rigorous error bound estimate, is called reduced basis method, as proposed by Y. Maday and E. M Rønquist (2002). The proper generalized decomposition (PGD) based model reduction, as an extension of POD, can be dated back to 1980s, and it was introduced by P. Ladevèze (1985) under the name of radial time-space approximation. F. Chinesta and his coworkers (2005, 2013) developed a PGD method to account for the parameter space, aiming at building offline computational vademecum for fast online predictions. It is noted that PGD methods are based on the idea of separation of variables and in particular a canonical tensor decompostion (CD). Zhang et al. (2021) consolidated the various attributes of the PGD methods with HiDeNN. The proposed HiDeNN-PGD method and the general reduced order machine learning finite element framework (Lu et al. 2021) provide powful tools for sovling large scale high dimensional problems.

The development of model reduction methods meets the urgent demand in the industry for fast and nearly real time simulations of engineering problems, such as online dynamic system control, structural health monitoring, vehicle health monitoring, on-line advanced manufacturing feedback control, automated driving controls and decisions, etc. Such applications usually require an intensive interaction between sensors, control algorithms, and simulation tools. Practical optimal control may require a reliable prediction within the range of milli- or sub-milli-second. Reducing the computation cost of simulations has been one of the major motivations for developing model reduction methods. Other reduced order modeling related topics are the feature engineering and data analytics, which constitute an extensive literature in field of machine learning. Thus, the reduced order modeling and the machine learning have intrinsic connections. Developing reduced order machine learning methods may enable physics-data combined models that can overcome the current bottleneck in model reduction methods and purely data-driven machine learning approaches.

### Acknowledgements
The authors would like to acknowledge the sources of some historic facts cited in this paper, which are taken from the following finite element history articles and books:

1. Oden, J.T. (1990). Historical comments on finite elements. In *A history of scientific computing* , 152-166.
2. Lax, P. (1993). Feng Kang. *SIAM News*, *26*(11).
3. Clough, R.W. and Wilson, E.L. (1999), August. Early finite element research at Berkeley. In *Fifth US National Conference on Computational Mechanics*, 1-35.
4. Felippa, C.A. (2004). Introduction to finite element methods. *University of Colorado*, *885*.
5. Owen, D.R.J. and Feng, Y.T. (2012). Fifty years of finite elements—a solid mechanics perspective. *Theoretical and Applied Mechanics Letters*, *2*(5), p.051001.
6. Stein, E. (2014). History of the finite element method–mathematics meets mechanics–part I: Engineering developments. In *The History of Theoretical, Material and Computational*